\documentclass[12pt]{article}
\usepackage[margin=2.5cm]{geometry}
\usepackage{amsmath}
\usepackage{amsfonts}
\usepackage{amsthm}
\usepackage{amssymb}
\usepackage{mathrsfs}
\usepackage{caption}
\usepackage{latexsym}
\usepackage{color}
\usepackage{hyperref}
\usepackage{url}
\usepackage{enumerate}
\usepackage{multicol} %����ʵ����ͬһҳ��ʵ�ֲ�ͬ�ķ���

\usepackage{multirow}
\usepackage{tikz}
\usepackage{bm}
\usepackage{booktabs}

\renewcommand\a{\alpha}    

\newcommand\Ga{\mathrm{\Gamma}}

\renewcommand\P{\mathrm{P}}     

\newcommand\ZZ{\mathbb{Z}}

          \newcommand\Aut{\mathrm{Aut}}
   \newcommand\Cay{\mathrm{Cay}}

\newcommand\soc{\mathrm{soc}}    \newcommand\Sy{\mathrm{S}}         
          
\newcommand\diam{\mathrm{diam}}            \newcommand\la{\langle}
\newcommand\ra{\rangle}

\newcommand\PGL{\mathrm{PGL}}                  \newcommand\PGammaL{\mathrm{P\Gamma L}}
                  
\newcommand\PSL{\mathrm{PSL}}

\newcommand\AS{\mathrm{AS}}

\newtheorem{theorem}{Theorem}[section]%
\newtheorem{theoremAB}{Theorem}

\newtheorem{lemma}[theorem]{Lemma}%
\newtheorem{proposition}[theorem]{Proposition}%
\begin{document}

\title{On the classification of $2$-distance-primitive graphs of small valency}

\author{
\parbox{\textwidth}{
\centering
Xing Zhang$^{\rm a}$, Xue Wang$^{\rm b}$\textcolor{blue}{,} Jun-Jie Huang$^{\rm c,*}$ \\[2pt]
\small
\hangindent=2em\hangafter=1
$^{\rm a}$ School of Mathematical Sciences, Zhejiang University, Hangzhou, 310058, P.R. China\\
$^{\rm b}$  School of Mathematical Science, Tianjin University of Technology, Tianjin, 300384, P.R. China\\
$^{\rm c}$School of Mathematics and Statistics, Chongqing University of Posts and Telecommunications, Chongqing, 400065, P. R. China }}

\renewcommand{\thefootnote}{\fnsymbol{footnote}}
\footnotetext[1]{Corresponding author. \\
E-mails: 0626554@zju.edu.cn(X.Zhang), xwang@email.tjut.edu.cn(X. Wang), huangjj@cqupt.edu.cn(J.-J. Huang)
}

\date{}
\maketitle

\begin{abstract}
A non-complete vertex-transitive graph is called
$2$-distance-primitive if a vertex stabilizer acts primitively on
both the first and second distance neighborhoods of the vertex.
Recently, Ruan, Liu and Jin [Discrete Math. 347 (2024), 113978]
showed that every connected $2$-distance-primitive graph of
valency at most $9$ either belongs to certain known families or
satisfies a collection of further restrictions, with the socle of
a vertex stabilizer isomorphic to $\PSL(2,7)$ or $\PSL(2,8)$.
In this paper, we resolve these remaining cases and thereby obtain
a complete classification of connected $2$-distance-primitive
graphs of valency at most $9$.
We also give a short proof of the main result of Ruan, Liu and Jin
and identify additional graphs omitted from their list.
\end{abstract}


\noindent\textbf{2020 MR Subject Classification:} 20B15, 20B30, 05C25.

\section{Introduction}\label{sec:introduction}
Throughout this paper, graphs are finite, connected, simple and
undirected. For a graph $\Gamma$, denote by $V(\Gamma)$ and
$\Aut(\Gamma)$ its vertex set and full automorphism group, respectively.
For $u,v\in V(\Gamma)$, let $d_\Gamma(u,v)$ denote the {\em distance} between
$u$ and $v$, and, for a nonnegative integer $i$, let
$\Gamma_i(u)=\{v\in V(\Gamma)\mid d_\Gamma(u,v)=i\}$.
In particular, we write $\Gamma(u)=\Gamma_1(u)$.
The {\em diameter} $\diam(\Ga)$ of a graph $\Ga$ is
the maximum distance between $u$ and $v$ for all $u,v\in V(\Ga)$.
For a finite group $G$,
denote by $\soc(G)$ the socle of $G$, that is, the subgroup
generated by all minimal normal subgroups of $G$.
Let $G$ be a transitive permutation group on the set $\Omega$.
Then $G$ is {\em primitive} if it has only trivial blocks in $\Omega$
(a {\em block} is a non-empty subset $\Delta$ of $\Omega$
such that $\Delta^g\cap\Delta\neq \emptyset$ implies $\Delta^g=\Delta$
for every $g\in G$ and a block $\Delta$ is {\em trivial} if $|\Delta| = 1$ or $|\Delta| = |\Omega|$).
It is well known that a transitive permutation group $G$ is primitive
if and only if $G_\a$ is a maximal subgroup of $G$, where $\a\in\Omega$;
see \cite[Corollary 1.5A]{Dixon}.

Let $G\leq\Aut(\Gamma)$ and let $1\leq s\leq\diam(\Gamma)$.
The graph $\Gamma$ is said to be \emph{$(G,s)$-distance-transitive}
if $G$ is transitive on $V(\Gamma)$ and, for each $u\in V(\Gamma)$,
the stabilizer $G_u$ is transitive on $\Gamma_i(u)$ for every
$1\leq i\leq s$.
If $s=\diam(\Gamma)$, then
$\Gamma$ is called \emph{$G$-distance-transitive}.
We omit $G$ from these terms when $G=\Aut(\Gamma)$.
Under this definition, a complete graph, which has diameter $1$,
is not $2$-distance-transitive. Every distance-transitive graph of diameter
at least $2$ is $2$-distance-transitive, but the converse does not hold.
A family of $2$-distance-transitive graphs that are not
distance-transitive was constructed in~\cite{HFZ}.
More generally, $s$-distance-transitive graphs have attracted considerable attention.
For example, Devillers, Giudici,
Li and Praeger~\cite{DGLP2012} developed a normal-quotient approach
to locally $s$-distance-transitive graphs and established reductions
to basic graphs with quasiprimitive or biquasiprimitive actions.
Moreover,
classifications of $2$-distance-transitive Cayley graphs
on cyclic and dihedral groups were obtained in~\cite{CJL,HFZY2025},
and $2$-distance-transitive graphs of valency at most $6$ were studied
in~\cite{CJS2017,JT16}.
For more results, refer to~\cite{DGLP13,HFZ} and the references therein.

A non-complete vertex-transitive graph $\Ga$ is called
$2$-distance-primitive if, for each $u\in V(\Ga)$, the stabilizer $A_u$,
where $A=\Aut(\Ga)$, is primitive on both $\Ga(u)$ and $\Ga_2(u)$.
This class of graphs was introduced by Jin et al.~\cite{JHL2019,JWZ2020}.
Clearly, every $2$-distance-primitive graph is $2$-distance-transitive.
It was proved in~\cite{JWZ2020} that if a $2$-distance-primitive graph
$\Ga$ is neither a cycle nor a complete bipartite graph, then $\Ga$ has
girth at most $4$, and the vertex stabilizer $A_u$ acts faithfully on both
$\Ga(u)$ and $\Ga_2(u)$. Building on this result,
$2$-distance-primitive graphs of prime valency were studied
in~\cite{HFKZY2024,JHL2019}. More recently, Ruan, Liu and
Jin~\cite{RLJ2024} investigated $2$-distance-primitive graphs of valency
at most $9$ and obtained the following result.

\begin{theoremAB}[{\cite[Theorem 1.2]{RLJ2024}}]\label{thm:A}
Let $\Gamma$ be a $2$-distance-primitive graph of
valency $k\leq9$. Let $A=\Aut(\Gamma)$ and $u\in V(\Gamma)$.
Then either $\Gamma$ is listed in Table~{\rm\ref{tab:graph}},
or $(k,c_2) = (8,2)$ and $\soc(A_u) \cong \PSL(2,7)$,
or $(k,c_2) = (9,2)$ and $\soc(A_u) \cong \PSL(2,8)$.
\end{theoremAB}

\begin{table}
\centering
\caption{$2$-distance-primitive graphs listed in \cite[Table~1]{RLJ2024}}\label{tab:graph}
\begin{tabular}{ll}\hline
valency $k$ & $\Gamma$ \\
\hline
2 & $C_n(n \ge 4)$ \\
3 & $K_{3,3}, K_{4,4} - 4K_2$ \\
4 & $K_{4,4}, K_{5,5} - 5K_2$ \\
5 & $K_{5,5}, K_{6,6} - 6K_2, H(5, 2), B(H_{11}), \text{icosahedron}$ \\
6 & $K_{6,6}, K_{7,7} - 7K_2, H(6, 2), B'(H_{11})$ \\
7 & $K_{7,7}, K_{8,8} - 8K_2, H(7, 2)$ \\
8 & $K_{8,8}, K_{9,9} - 9K_2, H(8, 2), B'(PG(3, 2))$ \\
9 & $K_{9,9}, K_{10,10} - 10K_2, H(9, 2), J(6, 3)$ \\
\hline
\end{tabular}
\end{table}

In this paper, we resolve the remaining cases in
Theorem~\ref{thm:A}. Our first theorem shows that no graph satisfying
either of the two exceptional conditions can occur.

\begin{theorem}\label{thm:exceptional}
There is no $2$-distance-primitive graph $\Gamma$ of valency $k$ satisfying
\[
(k,c_2,\soc(A_u))=(8,2,\PSL(2,7))
\quad\text{or}\quad (9,2,\PSL(2,8)),
\]
where $A=\Aut(\Gamma)$ and $u\in V(\Gamma)$.
\end{theorem}

Ruan et al.~\cite{RLJ2024} attempted to classify
$2$-distance-primitive graphs of valency at most $9$. However, their
proof contains several gaps. In particular, the folded $k$-cube
$\square_k$ (see~\cite[Lemma 2.10]{HFKZY2024}), where $k\geq5$, is a $2$-distance-primitive graph of
valency $k$ that was omitted from~\cite[Theorem~1.2]{RLJ2024}.
Moreover, the graph $B'(PG(3,2))$ listed in Table~\ref{tab:graph} is
not $2$-distance-primitive; see Lemma~\ref{exam:B32}. Accordingly, we
give a new proof of their classification theorem and obtain a corrected
complete classification of connected $2$-distance-primitive graphs of
valency at most $9$. Our proof of Theorem~\ref{thm:main} differs
substantially from that of~\cite[Theorem~1.1]{RLJ2024}.

\begin{theorem}\label{thm:main}
A non-complete graph $\Gamma$ of valency at most $9$ is
$2$-distance-primitive if and only if it is isomorphic to
one of the following:
\begin{enumerate}[\rm (1)]
\item $C_n$, where $n\geq4$;
\item $K_{k,k}$ or $K_{k+1,k+1}-(k+1)K_2$, where $3\leq k\leq9$;
\item $H(k,2)$ or $\square_k$, where $5\leq k\leq9$;
\item the icosahedron of valency $5$;
\item $B(H_{11})$ or $B'(H_{11})$, of valencies five and six, respectively;
\item $J(6,3)$, of valency nine.
\end{enumerate}
\end{theorem}

We remark that $H(3,2)\cong K_{4,4}-4K_2$ and $\square_4\cong K_{4,4}$
are already included, whereas $H(4,2)$ is not $2$-distance-primitive; see Proposition~\ref{examples}.
For descriptions of the graphs appearing in
Theorem~\ref{thm:main}, we refer the reader to Section~\ref{sec:graph}.

\section{Preliminaries and examples}\label{sec:prelim}

We use $\mathbb{Z}_n$ and $D_{2n}$ to denote the cyclic group
of order $n$ and the dihedral group of order $2n$, respectively.
For groups $N$ and $H$, write $N:H$ for a semidirect product
of $N$ by $H$.
For a permutation group $G$ on a set $\Omega$,
we denote by $G_{\a_1\a_2\cdots\a_k}$ the subgroup of $G$ that fixes $\{\a_1,\a_2,\cdots,\a_k\}$ pointwise,
where $\a_1,\a_2,\cdots,\a_k\in\Omega$.

Let $G$ be a primitive permutation group on a set $\Omega$ with socle $N$.
It is well known that finite primitive permutation groups fall into
eight types, of which the following two will be used in this paper.
The group $G$ is said to be {\em affine} (of type $\mathrm{HA}$)
if $N$ is abelian (and hence elementary abelian) and acts regularly
on $\Omega$, and {\em almost simple} (of type $\mathrm{AS}$)
if $N$ is a nonabelian simple group.
We refer the reader to~\cite{Praeger93,Praeger97} for the O'Nan--Scott
classification theorem and for further details.

\subsection{Graph-theoretic results}

For a positive integer $n$,
we denote by $C_n$ the cycle graph of length $n\geq 3$,
by $K_n$ the complete graph of order $n$,
by $K_{n,n}$ the complete bipartite graph of order $2n$,
and by $K_{n,n}-nK_2$ the graph obtained from
$K_{n,n}$ by deleting a perfect matching.
For a graph $\Ga$ with a subset $W$ of $V(\Gamma)$, denote by $\Gamma[W]$ the
subgraph of $\Gamma$ induced by $W$.

Let $s\geq1$ be an integer. An \emph{$s$-arc} of a graph $\Gamma$ is a sequence of $s+1$ vertices $(u_0,u_1,\ldots,u_s)$ in $\Gamma$ such that $u_i$ is adjacent to $u_{i+1}$ for $0\leq i\leq s-1$ and $u_{j-1}\neq u_{j+1}$ for $1\leq j\leq s-1$.
Let $G\leq\Aut(\Gamma)$. The graph $\Gamma$ is said to be \emph{$(G,s)$-arc-transitive} if $G$ is transitive on the set of all $s$-arcs of $\Gamma$.
Recall that a transitive permutation group $G$ on a set $\Omega$ is called
{\em $2$-transitive} if it acts transitively on the set of ordered pairs of distinct elements of $\Omega$.
If $\Gamma$ has valency at least two and $G$ is vertex-transitive, then $\Gamma$ is $(G,2)$-arc-transitive if and only if $G_u$ is $2$-transitive on $\Gamma(u)$ for some vertex $u$.
An $s$-arc $(u_0,u_1,\ldots,u_s)$ is called an
\emph{$s$-geodesic} if $d_\Gamma(u_0,u_s)=s$.
For $1\leq s\leq\diam(\Gamma)$, the graph $\Gamma$ is said to be
\emph{$(G,s)$-geodesic-transitive} if $G$ is transitive
on the set of all $i$-geodesics for each $1\leq i\leq s$.
When $G=\Aut(\Gamma)$, we simply say that $\Gamma$ is \emph{$s$-arc-transitive} and \emph{$s$-geodesic-transitive}, respectively.
Every $s$-geodesic-transitive graph is $s$-distance-transitive, whereas the converse does not hold; see~\cite{JDLP,Huang2027} for counterexamples.

\medskip
We begin with a structural result from
\cite[Theorem~1]{JWZ2020} that will be used repeatedly.

\begin{proposition}\label{faith:dis12}
Let $\Gamma$ be a $2$-distance-primitive graph. Then
either $\Gamma\cong C_n$ for some $n\geq4$, or $\Gamma$ is a complete
bipartite graph, or $\Gamma$ has girth at most $4$ and each vertex
stabilizer acts faithfully on both the first and second distance
neighbourhoods.
\end{proposition}

We next introduce the intersection parameters needed
in the subsequent counting arguments. Let $\Gamma$ be a
$(G,s)$-distance-transitive graph with $1\leq s\leq\diam(\Ga)$, and
let $u\in V(\Ga)$. For $1\leq i\leq s$ and $v\in\Ga_i(u)$, the
numbers of neighbours of $v$ in $\Ga_{i-1}(u)$, $\Ga_i(u)$, and
$\Ga_{i+1}(u)$ are independent of the choice of $v$; denote them by
$c_i$, $a_i$, and $b_i$, respectively. Thus $a_i+b_i+c_i$ is the
valency of $\Ga$. For a distance-transitive graph $\Ga$ of diameter
$d$, the sequence
$\{b_0,b_1,\ldots,b_{d-1};1,c_2,\ldots,c_d\}$ is called the
\emph{intersection array} of $\Ga$, where $b_0$ is the valency of
$\Ga$. These parameters play an important role in the study of
$s$-distance-transitive graphs; see, for example,
\cite{BCN1989,DGLP2012,JHL2019,JWZ2020}.

When $c_2$ is close to the valency, elementary
combinatorial arguments often determine the graph. The following
standard result records the two extremal cases; see also
\cite[Lemma~4.4]{CJS2017}.

\begin{proposition}\label{bipartgraph}
Let $\Gamma$ be a $2$-distance-transitive graph of valency $k\geq3$.
If $c_2=k$, then $\Gamma\cong K_{k,k}$; if $c_2=k-1$,
then $\Gamma\cong K_{k+1,k+1}-(k+1)K_2$.
\end{proposition}

The following elementary lemma collects two further
facts about intersection parameters that will be used in both main proofs.

\begin{lemma}\label{lem:basic}
Let $\Gamma$ be a $(G,s)$-distance-transitive graph of valency
$k$ and diameter $d$, where $2\leq s\leq d$. Then the following
hold.
\begin{enumerate}[\rm (i)]
\item For $u\in V(\Gamma)$ and $0\leq i<s$, $|\Gamma_i(u)|b_i=|\Gamma_{i+1}(u)|c_{i+1}$. Moreover,
$1=c_1\leq c_2\leq\cdots\leq c_s$, $k=b_0\geq b_1\geq\cdots\geq b_s$.
\item If $b_s\leq1$, then $\Gamma$ is $G$-distance-transitive.
\end{enumerate}
\end{lemma}

\proof Part (i) follows from~\cite[Lemma 2.4]{HFZ2026}.
To prove part (ii), let $u\in V(\Gamma)$. We first show that, for
$s\leq i<d$, each $v\in\Gamma_i(u)$ has at most $b_s$
neighbours in $\Gamma_{i+1}(u)$. Choose a vertex $w$ on a
geodesic from $u$ to $v$ such that $d_\Gamma(w,v)=s$, so that $d_\Gamma(u,w)=i-s$.
For $y\in\Gamma(v)\cap\Gamma_{i+1}(u)$, we have $s+1=d_\Gamma(u,y)-d_\Gamma(u,w)\leq d_\Gamma(w,y)\leq d_\Gamma(w,v)+1=s+1$.
Thus $\Gamma(v)\cap\Gamma_{i+1}(u)\subseteq \Gamma(v)\cap\Gamma_{s+1}(w)$, and the latter set has size $b_s\leq1$.

We now prove by induction on $i$ that $G_u$ is transitive on
$\Gamma_i(u)$ for $s\leq i\leq d$. The case $i=s$ follows
from the hypothesis. Suppose that $s\leq i<d$ and that
$G_u$ is transitive on $\Gamma_i(u)$. For
$y_1,y_2\in\Gamma_{i+1}(u)$, choose
$v_j\in\Gamma(y_j)\cap\Gamma_i(u)$ for $j=1,2$, and take
$g\in G_u$ with $v_1^g=v_2$. Both $y_1^g$ and $y_2$ are
neighbours of $v_2$ in $\Gamma_{i+1}(u)$, so the preceding
observation gives $y_1^g=y_2$. This completes the induction,
and hence $\Gamma$ is $G$-distance-transitive. \qed

\medskip
A regular graph on $v$ vertices with valency $k$ is
called \emph{amply regular} with parameters $(v,k,\lambda,\mu)$ if
any two adjacent vertices have precisely $\lambda$ common neighbours
and any two vertices at distance $2$ have precisely $\mu$ common
neighbours. Thus every $2$-distance-transitive graph is amply regular
with $\lambda=a_1$ and $\mu=c_2$. Standard properties of amply regular graphs
(see~\cite[Theorem~1.5.5 and Proposition~1.9.1]{BCN1989})
yield the following result.

\begin{proposition}\label{amply}
Let $\Gamma$ be a connected $2$-distance-transitive graph and let $u\in V(\Ga)$.
\begin{enumerate}[\rm (1)]
  \item $|\Ga(w)\cap\Ga_2(u)|\geq 3$ for every $w\in\Ga_3(u)$.
  \item If $\Ga$ has diameter at least $3$, then $|\Ga_2(u)|\geq |\Ga(u)|$ and this equality holds if and only if $\Gamma$ is a polygon or a Taylor graph.
\end{enumerate}
\end{proposition}

Here a polygon is simply a cycle $C_n$ with $n\geq5$. A \emph{Taylor graph} is a distance-regular graph of diameter $3$ with intersection array $\{k,\mu,1;1,\mu,k\}$; see~\cite[p.~13]{BCN1989}.

\subsection{Examples of $2$-distance-primitive graphs}\label{sec:graph}
Let $d,n\geq2$ be integers.
The {\em Hamming graph} $H(d,n)$ is the graph with vertex set
$\{(x_1,x_2,\ldots,x_d)\mid x_i\in\{0,1,\ldots,n-1\}
\text{ for }1\leq i\leq d\}$, where two vertices $(x_1,x_2,\ldots,x_d)$ and $(y_1,y_2,\ldots,y_d)$ are adjacent if and only if $|\{i\mid x_i\neq y_i,1\leq i\leq d\}|=1$.
In particular, $H(d,2)$ is known as the {\em $d$-cube},
denoted by $Q_d$ in some references, and $H(2,2)\cong C_4$.

Let $d\geq3$ be an integer.
The {\em folded $d$-cube} $\square_d$ is the graph obtained from
the $d$-cube $H(d,2)$ by identifying antipodal pairs of vertices,
that is, its vertex set is
$\{\{u,u'\}\mid u,u'\in V(H(d,2)),\ d_{H(d,2)}(u,u')=d\}$.
Two vertices $\{u,u'\}$ and $\{v,v'\}$ are adjacent in $\square_d$
if and only if there exist $x\in\{u,u'\}$ and $y\in\{v,v'\}$
such that $x$ and $y$ are adjacent in $H(d,2)$.
Note that $H(d,2)$ is a double cover of $\square_d$.
In particular, $\square_3\cong K_4$ and $\square_4\cong K_{4,4}$.
The $2$-distance-primitivity of $\square_d$ was determined
in~\cite[Lemma~2.10]{HFKZY2024}.

We next turn to Hamming graphs and record several standard examples, as summarised in~\cite[Section~2.3]{HFKZY2024}.

\begin{proposition}\label{examples}
The following graphs are $2$-distance-primitive: $C_n$ with $n\geq 4$, $K_{n,n}$ with $n\geq 3$,
$K_{n,n}-n K_2$ with $n\geq 4$, the Hamming graph $H(d,2)$ with $d\geq2$ and $d\neq4$,
the folded $d$-cube $\square_d$ with $d\geq 4$, and the icosahedron.
\end{proposition}

We remark that the icosahedron is denoted by $\mathrm{Cos}(5)$ in~\cite[Definition~2.5]{HFKZY2024}.
The $2$-distance-primitivity of Hamming graphs $H(d,2)$ was established in~\cite[Proposition~3.2]{RLJ2024},
and we provide a simpler proof that differs from the one in~\cite{RLJ2024}.

\begin{lemma}\label{exam:Hamming}
The Hamming graph $H(d,n)$ with $d,n\geq2$ is $2$-distance-primitive if and only if $n=2$ and $d\neq4$.
\end{lemma}

\proof Let $\Ga= H(d,n)$,
and let $A=\Aut(\Gamma)$.
By~\cite[Theorem~9.2.1]{BCN1989},
we have $A\cong S_n\wr S_d$, and $\Gamma$ is distance-transitive.
Let $u\in V(\Ga)$.
Suppose that $n\geq3$.
Then $\Gamma[\Gamma(u)]\cong dK_{n-1}$, whose connected
components form an $A_u$-invariant partition of $\Gamma(u)$.
Since $d,n-1\geq2$, this partition is nontrivial. Thus
$A_u$ is imprimitive on $\Gamma(u)$, and $\Gamma$ is not
$2$-distance-primitive.
Therefore $n=2$, and the result follows from
Proposition~\ref{examples}.
\qed

\medskip
We next turn to Johnson graphs.
Let $\Omega$ be a finite set of cardinality $n$,
and let $\Omega^{\{k\}}$ be the set of all $k$-subsets of $\Omega$.
The {\em Johnson graph $J(n,k)$} is the graph whose vertex set is $\Omega^{\{k\}}$,
and two $k$-subsets $U$ and $W$ are adjacent if and only if $|U\cap W|=k-1$.
It is well known that $J(n,k)\cong J(n,n-k)$ and that $J(n,1)\cong K_n$.
Moreover, the Johnson graphs $J(n,2)$ are also called \emph{triangular graphs}.
The $2$-distance-primitivity of Johnson graphs was studied in~\cite{RLJ2024}.
For completeness, we also give an alternative proof that differs from the one in~\cite{RLJ2024}.

\begin{lemma}\label{exam:Johnson}
Let $2\leq k\leq n/2$. Then the Johnson graph $J(n,k)$
is $2$-distance-primitive if and only if $n=2k$ and either
$k=3$ or $k\geq5$.
\end{lemma}

\proof
For a positive integer $n$, let $[n]=\{1,2,\ldots,n\}$.
Let $\Omega=[n]$ and let $\Gamma=J(n,k)$.
Let $u\in \Omega^{\{k\}}$ and let $\overline{u}=\Omega\setminus u$.
For $i=1,2$, every vertex in
$\Gamma_i(u)$ can be written uniquely in the form
\[
 (u\setminus X)\cup Y,
 \quad \text{with~}
 X\in u^{\{i\}} \text{~and~}
 Y\in\overline{u}^{\{i\}}.
\]
Consequently, $\Gamma_i(u)$ can be identified with $u^{\{i\}}\times \overline{u}^{\{i\}}$.
Let $A=\Aut(\Gamma)$.
By~\cite[Theorem~9.1.2]{BCN1989},
\[
 \Aut(J(n,k))\cong
 \begin{cases}
  S_n, & n\neq 2k,\\
  S_n\times \mathbb{Z}_2, & n=2k,
 \end{cases}
\]
where $S_n$ acts on $\Omega^{\{k\}}$ via its natural action on
$\Omega$, and in the second case the additional involution is the
complementation map $u\mapsto \overline{u}$.

Suppose first that $n>2k$. Then $A_u\cong S_k\times S_{n-k}$
and, under the identification
$\Gamma(u)= u\times\overline{u}$, the group $A_u$ acts
coordinatewise.
Therefore, $\bigl\{\{x\}\times\overline{u}\mid x\in u\bigr\}$
is a nontrivial $A_u$-invariant partition of $\Gamma(u)$ into $k$ blocks of size $n-k$.
Thus $A_u$ is imprimitive on $\Gamma(u)$, and so $\Gamma$ is
not $2$-distance-primitive. Therefore, it is necessary that $n=2k$.

Now suppose that $n=2k$. Then $A_u\cong S_k\wr S_2$, where an element of $A_u$ either stabilises $u$ setwise, or interchanges
$u$ and $\overline{u}$.
Under the identification $\Ga_i(u)\cong u^{\{i\}}\times\overline{u}^{\{i\}}$,
the action of $A_u$ on $\Gamma_i(u)$ is permutation-isomorphic to
the product action of $S_k\wr S_2$ on $[k]^{\{i\}}\times[k]^{\{i\}}$.
Recall that if $H$ is a transitive group on a set $\Delta$,
then $H\wr S_2$ in its product action on $\Delta^2$ is primitive
if and only if $H$ is primitive and nonregular on $\Delta$;
see~\cite[Lemma~2.7A]{Dixon}.

For $i=1$, the natural action of $S_k$ on $[k]$ is primitive,
and it is nonregular precisely when $k\geq3$. Hence $A_u$ is
primitive on $\Gamma(u)$ if and only if $k\geq3$. Consequently, $J(4,2)$ is not $2$-distance-primitive.

It remains to consider the action of $A_u$ on $\Gamma_2(u)$,
identified with $[k]^{\{2\}}\times[k]^{\{2\}}$.
The action of $A_u$ is then the product action arising from the
action of $S_k$ on $[k]^{\{2\}}$.
If $k=3$, then the action of $S_3$ on its $2$-subsets is
permutation-isomorphic to its natural action of degree $3$
and is therefore primitive and nonregular.
If $k\geq5$, the stabiliser in $S_k$ of a $2$-subset is
$S_2\times S_{k-2}$, which is a maximal subgroup of $S_k$
since $2<k/2$; see~\cite[Exercise~5.2.2]{Dixon}.
Thus the action of $S_k$ on $[k]^{\{2\}}$ is primitive,
and it is clearly nonregular. Hence $A_u$ is primitive on
$\Gamma_2(u)$.

Finally, suppose that $k=4$.
The partition $ \left\{\{X,[4]\setminus X\}\mid X\in[4]^{\{2\}}\right\}$
is a nontrivial $S_4$-invariant partition of $[4]^{\{2\}}$ into three complementary pairs.
Therefore, $S_4$ is imprimitive on its $2$-subsets, and the corresponding
product action of $S_4\wr S_2$ on $\Gamma_2(u)$ is also
imprimitive. Hence $k=4$ must be excluded.
We therefore conclude that $A_u$ is primitive on both $\Gamma(u)$ and
$\Gamma_2(u)$ if and only if
\[
 n=2k
 \quad\text{and}\quad
 k=3\ \text{or}\ k\geq5.
\]
Therefore, $J(n,k)$ is $2$-distance-primitive precisely in these cases. \qed

\medskip
We next recall the Cayley graph description of two sporadic examples used below.
Let $G$ be a finite group and let
$S\subseteq G\setminus\{1\}$ satisfy $\la S\ra=G$ and $S=S^{-1}$.
The {\em Cayley graph} $\Cay(G,S)$ on $G$ with respect to $S$ is
the graph with vertex set $G$ and edge set
$\{\{g,sg\}\mid g\in G,s\in S\}$.

\medskip
Let $B(H_{11})$ and $B'(H_{11})$ denote, respectively,
the incidence and nonincidence graphs of the unique
$2$-$(11,5,2)$ design. These are distance-transitive graphs of
valencies $5$ and $6$, respectively.
Let
\begin{eqnarray*}
&G=\langle a,b\mid a^{11}=b^2=1,\ bab=a^{-1}\rangle\cong D_{22},\\
&S_1=\{ab,a^3b,a^4b,a^5b,a^9b\},\,S_2=\{b,a^2b,a^6b,a^7b,a^8b,a^{10}b\}\textcolor{blue}{.}
\end{eqnarray*}
By Magma~\cite{Magma}, $B(H_{11})\cong\Cay(G,S_1)$ and
$B'(H_{11})\cong\Cay(G,S_2)$.

\begin{lemma}\label{exam:B11}
Both $B(H_{11})$ and $B'(H_{11})$ are $2$-distance-primitive.
\end{lemma}

\proof That $B(H_{11})$ is $2$-distance-primitive
follows from~\cite[Example~2.7]{HFKZY2024}.
We now consider the $2$-distance-primitivity of $B'(H_{11})$.
Let $\Gamma=B'(H_{11})$ and let $A=\Aut(\Ga)$.
For $u\in V(\Ga)$, we have
$A\cong\PGL(2,11)$ and $A_u\cong A_5$. Moreover, $\Ga$ is
distance-transitive with intersection array $\{6,5,3;1,3,6\}$;
see~\cite[Table~1]{Huang2027}, where it is denoted by
$\mathcal{G}_{22,6}$.
By Lemma~\ref{lem:basic}(1), $|\Ga_2(u)|=10$.
Clearly, the action of $A_u$ on $\Ga(u)$ is $2$-transitive
(the unique transitive action of $A_5$ of degree $6$, with point
stabilizer $D_5$); in particular, $A_u$ is primitive on $\Ga(u)$.
Let $w\in\Ga_2(u)$. Since $A_u$ is transitive on $\Ga_2(u)$,
$|A_{uw}|=|A_u|/|\Ga_2(u)|=6$; as $A_5$ has no element of order $6$, $A_{uw}\cong S_3$.
As $S_3$ is maximal in $A_5$, the action of $A_u$ on $\Ga_2(u)$
is primitive.
Thus, $\Gamma$ is $2$-distance-primitive. \qed

\medskip
The graph $\Ga:=B'(PG(3,2))$ is the
point--hyperplane nonincidence graph of $PG(3,2)$. It is
distance-transitive with intersection array $\{8,7,4;1,4,8\}$;
see~\cite[Table~1]{Huang2027}, where it is denoted by
$\mathcal{G}_{30,8}$. Moreover,
$\Aut(\Ga)\cong\PSL(4,2).\mathbb{Z}_2$, and for $u\in V(\Ga)$,
the vertex stabilizer satisfies
$A_u\cong\mathbb{Z}_2^3:\PSL(3,2)$.
In particular, $\Ga$ is neither a cycle nor a complete bipartite graph.
If $\Ga$ is $2$-distance-primitive,
then by Proposition~\ref{faith:dis12},
we conclude that $A_u$ acts faithfully and primitively on $\Ga(u)$ and $\Ga_2(u)$.
Thus, the normal subgroup isomorphic to $\mathbb{Z}_2^3$ in $A_u$
must act regularly on $\Ga(u)$ and $\Ga_2(u)$.
This implies that $|\Ga(u)|=|\Ga_2(u)|=8$.
However, by Lemma~\ref{lem:basic}(1), $|\Ga_2(u)|=14$, leading to a contradiction.
This proves the following lemma.

\begin{lemma}\label{exam:B32}
The graph $B'(PG(3, 2))$ is not $2$-distance-primitive.
\end{lemma}

We conclude this section by recalling two notions
that will be used repeatedly.
A {\em rectagraph} is a connected triangle-free graph
with the property that every path of length two is in a unique quadrangle;
see~\cite{BCN1989}.
A {\em $(0,2)$-graph} is a connected graph such that
any two distinct vertices have either $0$ or $2$ common neighbours.
Clearly, a $2$-distance-transitive graph, with girth $4$ and $c_2=2$,
is a rectagraph and a $(0,2)$-graph.

\section{Proofs of Theorems~\ref{thm:exceptional} and \ref{thm:main}}\label{sec:exceptional}

\medskip\noindent{\bf Proof of Theorem~\ref{thm:exceptional}\textcolor{blue}{.}}
Let $\Gamma$ be a $2$-distance-primitive graph of valency $k$, let
$A=\Aut(\Gamma)$, and fix $u\in V(\Gamma)$.
Assume that
\begin{equation*}
(k,c_2,\soc(A_u))=(8,2,\PSL(2,7))\text{~or~} (9,2,\PSL(2,8)).
\end{equation*}
Then $\Gamma$ is neither a cycle nor a complete bipartite graph.
By Proposition~\ref{faith:dis12},
$A_u$ acts faithfully on both $\Gamma(u)$ and $\Gamma_2(u)$,
and these actions are primitive by the $2$-distance-primitivity of $\Gamma$.
In the two cases under consideration,
the action of $A_u$ on $\Gamma(u)$ is $2$-transitive.
Thus the induced graph $\Gamma[\Gamma(u)]$ is either complete or edgeless.
In the former case, $\{u\}\cup\Gamma(u)$ induces a complete graph,
and so the connectivity of $\Gamma$ would imply that $\Gamma\cong K_{k+1}$, a contradiction.
Thus $\Gamma[\Gamma(u)]$ is edgeless, and therefore $b_1=k-1$.
Since $c_2=2$, it follows that $\Gamma$ has girth $4$,
and every two distinct vertices have either $0$ or $2$ common neighbours.
Thus $\Gamma$ is a rectagraph and, in particular, a $(0,2)$-graph.
The $(0,2)$-graphs of valency at most $8$ were classified in~\cite{Bro,BO2009}.
Let $(u,v,w)$ be a $2$-geodesic of $\Gamma$.
We shall show that no graph satisfying the above
assumptions exists.

\medskip\noindent{\bf Case 1: $k=8$.}

In this case, $A_u\cong \PSL(2,7)$ or $\PGL(2,7)$.
By Lemma~\ref{lem:basic}(1), we have $|\Gamma_2(u)|=8b_1/c_2=28$.
Since $A_u$ acts primitively on $\Gamma_2(u)$,
the stabiliser $A_{uw}$ is a maximal subgroup of $A_u$ of index $28$
(see~\cite[Corollary~1.5A]{Dixon}).
However, $\PSL(2,7)$ has no maximal subgroup of index $28$
according to the Atlas~\cite{Atlas}.
Therefore, $A_u\cong\PGL(2,7)$.
As $\Gamma$ is a $(0,2)$-graph,
$\Gamma$ is one of the $302$ connected $(0,2)$-graphs classified in~\cite[Tables~1--2]{BO2009}.
In those tables, $v$ denotes the number of vertices, $\mathrm{gsz}$ denotes
the order of the full automorphism group,
and $\mathrm{tra}$ indicates vertex-transitivity.
Since $A$ is vertex-transitive, the orbit--stabiliser theorem yields
$\mathrm{gsz}/v=|A|/|V(\Gamma)|=|A_u|=336$.

Suppose first that $\Gamma$ is bipartite.
Then $\Gamma$ occurs among the $104$ graphs in~\cite[Table~1]{BO2009}.
Among the vertex-transitive entries, none satisfies $\mathrm{gsz}/v=336$, a contradiction.

Suppose next that $\Gamma$ is non-bipartite.
Then $\Gamma$ occurs among the $198$ graphs in~\cite[Table~2]{BO2009}.
Since $\Gamma$ is vertex-transitive and a rectagraph,
its entry must have $\mathrm{Orbits}=\mathrm{tra}$ and
$\mathrm{Locally}=\mathrm{recta}$.
Among the entries satisfying these conditions, the possible
values of $\mathrm{gsz}/v$ are $6$, $18$, $20$, $24$, $40$, $48$, $120$, $144$, $720$, $5040$.
None equals $336$, again a contradiction.

\medskip\noindent{\bf Case 2: $k=9$.}

In this case, $A_u\cong\PSL(2,8)$ or $\PGammaL(2,8)$.
Again by Lemma~\ref{lem:basic}(1),
we have $|\Gamma_2(u)|=kb_1/c_2=9\cdot8/2=36$.
If $\Gamma$ is distance-transitive,
then by \cite[Theorems~7.5.2 and~7.5.3(vi)]{BCN1989} or \cite[Table~1]{Huang2027},
the only distance-transitive graphs of valency $9$, girth $4$, and with $c_2=2$ are $H(9,2)$ and $\square_9$.
Both have full vertex stabiliser $S_9$, contrary to
$A_u\cong\PSL(2,8)$ or $\PGammaL(2,8)$.
Therefore, $\Gamma$ is not distance-transitive.

Assume that $a_2\neq 0$.
Then the subgraph $\Gamma[\Gamma_2(u)]$ is a
vertex-transitive graph of valency $a_2$,
and $\Gamma(w)\cap\Gamma_2(u)$ is a union of nontrivial orbits of $A_{uw}$.
If $A_u\cong\PSL(2,8)$,
then a computation in Magma~\cite{Magma} shows that $A_{uw}\cong D_{14}$
and $A_{uw}$ has one trivial orbit (namely, $\{w\}$) of length $1$,
three orbits of length $7$\textcolor{blue}{,} and one orbit of length $14$
on $\Gamma_2(u)$.
Since $a_2+b_2+c_2=k$ and $c_2=2$, we have $a_2\leq 7$, and so $a_2=7$.
It follows that $(a_2,b_2,c_2)=(7,0,2)$, and hence $\Gamma$ has diameter $2$.
Since $\Gamma$ is $2$-distance-primitive, it is
distance-transitive of order $1+9+36=46$, contrary to the
classification in~\cite[Theorems~7.5.2 and~7.5.3]{BCN1989}.
If $A_u\cong \PGammaL(2,8)$, then a computation in
Magma~\cite{Magma} shows that $A_{uw}\cong \mathbb{Z}_7:\mathbb{Z}_6$,
and $A_{uw}$ has two nontrivial orbits, of lengths
$14$ and $21$, on $\Gamma_2(u)$.
This implies that if $a_2\neq 0$, then $a_2\geq 14$, which is impossible.
We therefore conclude that $a_2=0$ and $b_2=7$.

Write $\Gamma(u)\cap\Gamma(w)=\{v,z\}$.
Let $P$ be a subgroup of order $7$ in $A_{uw}$.
Since $P$ preserves the $2$-set $\{v,z\}$, it fixes both $v$ and $z$,
and hence $P\leq A_{uvw}$.
The group $P$ preserves
$\Gamma(w)\cap\Gamma_3(u)=\Gamma(w)\setminus\{v,z\}$,
which has size $7$. If $P$ fixed this set pointwise, then it would fix
$\Gamma(w)$ pointwise, contradicting the faithfulness of the local
action at $w$. Thus $P$ is transitive (indeed, regular) on this set.
It follows that $\Gamma$ is $(A,3)$-geodesic-transitive; see
\cite[Corollary~2.5]{Huang2027}. Consequently, $\Gamma$ is also
$(A,3)$-distance-transitive.
By Lemma~\ref{lem:basic}(1), we have
$|\Gamma_3(u)|c_3=|\Gamma_2(u)|b_2=36\cdot7=252$.
In particular, $c_3\mid252$.
By Proposition~\ref{amply}~(1), we have $c_3\geq3$.
If $b_3\leq1$, then Lemma~\ref{lem:basic}(2) implies that $\Gamma$ is distance-transitive,
a contradiction.
Hence $b_3\geq2$. Since
$a_3+b_3+c_3=9$, it follows that $3\leq c_3\leq7$,
and therefore $c_3\in\{3,4,6,7\}$.
Let $x\in\Gamma(w)\cap\Gamma_3(u)$.

Suppose first that $c_3=7$. Then $|\Gamma_3(u)|=36$,
$a_3+b_3=2$, and hence $(a_3,b_3)=(0,2)$.
Put
\[
 \Delta_1:=\Gamma(x)\cap\Gamma_2(u),
 \quad
 \Delta_2:=\Gamma(x)\cap\Gamma_4(u).
\]
Then $|\Delta_1|=7$ and $|\Delta_2|=2$.
Note that $|A_{ux}|=14$ or $42$.
Let $P=\langle\alpha\rangle$ be a subgroup of order $7$ in $A_{ux}$.
Since $|\Delta_2|=2$, $P$ fixes $\Delta_2$ pointwise.
On the other hand,
Proposition~\ref{faith:dis12} implies that $A_x$ acts faithfully on $\Gamma(x)$.
Hence the action of $P$ on $\Gamma(x)$ is faithful,
so $P$ cannot fix $\Delta_1$ pointwise.
It follows that $P$ acts regularly on $\Delta_1$.

If $A_u\cong\PGammaL(2,8)$, then a computation using
Magma~\cite{Magma} shows that
$A_{ux}\cong\mathbb{Z}_7:\mathbb{Z}_6$ and that its orbits on
$\Gamma_2(u)$ have lengths $1$, $14$, and $21$.
This is impossible because $\Delta_1$ is an $A_{ux}$-invariant subset of
$\Gamma_2(u)$ of size $7$.
It remains to consider $A_u\cong\PSL(2,8)$.
Again by Magma~\cite{Magma}, $A_{ux}\cong D_{14}$.
In the natural action of $A_x\cong\PSL(2,8)$ on $\Gamma(x)$ of
degree $9$, every subgroup isomorphic to $D_{14}$ has two orbits,
of lengths $2$ and $7$. Thus $A_{ux}$ has two orbits of these lengths
on $\Gamma(x)$.
The regularity of $P$ on $\Delta_1$ shows that $A_{ux}$ is transitive on $\Delta_1$,
and hence $A_{ux}$ is also transitive on $\Delta_2$.
Since $\Gamma$ is $(A,3)$-distance-transitive,
$A_u$ is transitive on $\Gamma_3(u)$. Together with the transitivity
of $A_{ux}$ on $\Delta_2=\Gamma(x)\cap\Gamma_4(u)$, this implies
that $A_u$ is transitive on $\Gamma_4(u)$. Hence $\Gamma$ is
$(A,4)$-distance-transitive. By Lemma~\ref{lem:basic}(1),
we have $|\Gamma_4(u)|c_4=|\Gamma_3(u)|b_3=72$ and
$7=c_3\leq c_4\leq9$. Thus $c_4\in\{8,9\}$, so $b_4\leq1$.
Lemma~\ref{lem:basic}(2) now implies that $\Gamma$ is
distance-transitive, again a contradiction.

Suppose next that $c_3=6$. Then $|\Gamma_3(u)|=42$.
However, a computation in Magma~\cite{Magma} shows
that $A_u$ has no subgroup of index $42$, a contradiction.

Suppose that $c_3=4$. Then $|\Gamma_3(u)|=63$, and $|A_{ux}|=8$ or $24$.
Let $L$ be a Sylow $2$-subgroup of $A_{ux}$.
Since $|A_u:\soc(A_u)|$ is odd, $L\leq\soc(A_u)$ and $|L|=8$; thus
$L$ is a Sylow $2$-subgroup of $\PSL(2,8)$.
A computation in Magma~\cite{Magma} shows that, in
the natural action of $A_x$ on $\Gamma(x)$, every Sylow $2$-subgroup
of $A_x$ has two orbits, of lengths $1$ and $8$.
Consequently, the orbit lengths of $L$ on
$\Gamma(x)$ are $1$ and $8$.
This is impossible, since $L\leq A_{ux}$ stabilises
the $4$-set $\Gamma(x)\cap\Gamma_2(u)$ setwise.
Therefore $c_3\neq4$.

Having excluded $c_3\in\{4,6,7\}$, we must have $c_3=3$.
Then $\Gamma$ is a rectagraph with $a_2=0$ and $c_3=3$.
Moreover, $A_u$ is $4$-homogeneous on $\Gamma(u)$;
see~\cite[Theorem~1(iv)]{Kantor1972}.
By~\cite[Theorem~1.3]{BDFP2015},
$\Gamma\cong Q_9$ or $\square_9$, since the other graphs in that theorem
have valency $23$ or $24$.
Both possibilities have already been excluded.
This final contradiction completes the proof.\hfill\qed

\medskip
The preceding proof rules out the two exceptional cases left open
in Theorem~\ref{thm:A}. We can now complete the classification.
We reconsider the possible vertex stabilizers, taking into account
the corrections described in the introduction, and verify that every
$2$-distance-primitive graph of valency at most $9$ indeed appears
in the list in Theorem~\ref{thm:main}.

\medskip\noindent{\bf Proof of Theorem~\ref{thm:main}\textcolor{blue}{.}}
Let $\Gamma$ be a finite connected non-complete $2$-distance-primitive graph of valency $k\leq9$.
Let $A=\Aut(\Gamma)$ and let $u\in V(\Gamma)$.
If $\Gamma$ is a cycle or a complete bipartite graph,
then it appears in part~(1) or~(2), respectively.
In what follows, we always assume that $\Gamma$ is neither a cycle nor a complete bipartite graph.
By Proposition~\ref{faith:dis12}, $A_u$ acts
faithfully on $\Gamma(u)$ and $\Gamma_2(u)$; by hypothesis, both
actions are primitive.
By the classification of primitive groups of degree at most $9$
(see~\cite{Dixon}; this can also be checked using
Magma~\cite{Magma}), we conclude that $A_u$ is of type $\mathrm{HA}$
or $\mathrm{AS}$ and appears in Table~\ref{ASgroup}.

Put $T=\operatorname{soc}(A_u)$. In the affine case, $T$ is
elementary abelian and regular on $\Gamma(u)$. In the almost
simple case, $T$ is nonabelian simple and
$T\leq A_u\leq\Aut(T)$.
We treat the affine and almost simple cases separately.

\begin{table}
\centering
\caption{Almost simple primitive groups of degree at most $9$}\label{ASgroup}
\begin{tabular}{cl}
\toprule
Degree $k$ & $A_u$\\
\midrule
$5\leq k\leq9$ & $A_k$, $S_k$ in their natural actions\\
$6$ & $\PSL(2,5)$, $\PGL(2,5)$\\
$7$ & $\PSL(2,7)$\\
$8$ & $\PSL(2,7)$, $\PGL(2,7)$\\
$9$ & $\PSL(2,8)$, $\PGammaL(2,8)$\\
\bottomrule
\end{tabular}
\end{table}

\medskip
\noindent
{\bf Case 1.} $A_u$ is of type $\mathrm{HA}$.

Since $A_u$ is faithful and primitive on both $\Gamma(u)$ and $\Gamma_2(u)$,
the socle $T$ is transitive on each of $\Gamma(u)$ and $\Gamma_2(u)$,
and hence is regular on each set because $T$ is abelian.
Thus, $|\Gamma(u)|=|T|=|\Gamma_2(u)|=k$.
By Lemma~\ref{lem:basic}(1), we have $b_1=c_2$.

Assume that $\Gamma$ has girth $4$.
Then $b_1=c_2=k-1$.
Since $\Ga$ is $2$-distance-transitive,
it follows from Proposition~\ref{bipartgraph} that $\Gamma\cong K_{k+1,k+1}-(k+1)K_2$ with $k\geq3$.
Now, $A\cong\Sy_{k+1}\times\mathbb{Z}_2$ and $A_u\cong\Sy_k$.
Since $A_u\cong S_k$ is of type $\mathrm{HA}$ and
$k\geq3$, we have $k\in\{3,4\}$. Hence
$\Ga\cong K_{4,4}-4K_2$ or $K_{5,5}-5K_2$, both of which are
listed in part~(2).

Assume that $\Gamma$ has girth $3$.
If $\Gamma$ has diameter at least $3$,
then Proposition~\ref{amply}(2) implies that
$\Gamma$ is a Taylor graph of intersection array $\{k,c_2,1;1,c_2,k\}$.
In particular, $b_2=1$,
and so $\Gamma$ is distance-transitive by Lemma~\ref{lem:basic}(2).
Since $\Ga$ has girth $3$ and valency $k\leq9$,
the classifications in~\cite[Theorems~7.5.2 and~7.5.3]{BCN1989}
or~\cite[Table~1]{Huang2027} show that the only possibilities are
the icosahedron and $J(6,3)$, which appear in parts~(4) and~(6),
respectively.

Now suppose that $\Gamma$ has diameter $2$. Then $\Gamma$ is
distance-transitive with intersection array $\{k,b_1;1,b_1\}$.
By the same classifications, the graphs
of girth three and valency at most nine with this
intersection-array form are $H(2,3)\cong\P(9)$, $\P(13)$, and $\P(17)$.
If $\Ga\cong H(2,3)$, then by Lemma~\ref{exam:Hamming}, $\Gamma$ is not
$2$-transitive primitive, a contradiction. If $\Ga\cong\P(q)$ with $q\in\{13,17\}$, then
$A\cong\ZZ_q:\ZZ_{(q-1)/2}$ and $A_u\cong\ZZ_{(q-1)/2}$ acts
regularly on $\Ga(u)$. As $(q-1)/2\in\{6,8\}$ is composite,
this action is imprimitive, again a contradiction. Therefore the case $\diam(\Gamma)=2$ does not occur.

\medskip
\noindent
{\bf Case 2.} $A_u$ is of type $\mathrm{AS}$.

In this case, $A_u$ is listed in Table~\ref{ASgroup}, and is $2$-transitive on $\Gamma(u)$.
Since $\Gamma$ is non-complete, its local graph must therefore be empty.
Hence $\Gamma$ has girth $4$ and is $(A,2)$-arc-transitive.
In particular, $b_1=k-1$,
$|\Gamma_2(u)|c_2=k(k-1)$, and $2\leq c_2\leq k$.

The cases $c_2=k$ and $c_2=k-1$ are settled by
Proposition~\ref{bipartgraph}. They yield $\Ga\cong K_{k,k}$, which
has already been excluded by our standing assumption, or
$\Ga\cong K_{k+1,k+1}-(k+1)K_2$. In the latter case, $k\geq5$
because $A_u$ is of type $\AS$, and the graph appears in part~(2).
In what follows, assume that $2\leq c_2<k-1$.

Suppose that $A_u$ induces $A_k$ or $S_k$ in its natural
action on $\Gamma(u)$. By~\cite[Proposition~2.3]{HFKZY2024},
which records~\cite[Theorem~4.5]{Cameron1975},
$\Gamma$ is isomorphic to $H(k,2)$, $\square_k$, or
$B(H_{11})$ with $k=5$.
Such graphs are listed in parts (3) and (5).

For the remaining cases, suppose first that $k=7$.
Then $A_u\cong\PSL(3,2)\cong\PSL(2,7)$.
Since $c_2\mid42$ and $2\leq c_2<6$, we have $c_2=2$ or $3$,
and hence $|\Gamma_2(u)|=21$ or $14$, respectively.
However, by Magma~\cite{Magma} or Atlas~\cite{Atlas}, the maximal
subgroups of $\PSL(2,7)$ have indices $7$ and $8$.
Thus $A_u$ has no primitive permutation representation
of degree $14$ or $21$, a contradiction.

It remains to consider $k\in\{6,8,9\}$. Put
$q=k-1$; then $\PSL(2,q)\leq A_u\leq\PGammaL(2,q)$ in its natural
action on the projective line.
Suppose that $c_2=2$.
If $k=6$, then $A_u\cong A_5$ or $S_5$ and
$|\Gamma_2(u)|=15$.
The maximal subgroups of $A_5$ have indices $5,6,10$, while those
of $S_5$ have indices $2,5,6,10$.
Neither group admits a primitive action of degree $15$,
a contradiction.
If $k=8$ or $9$, the corresponding possibilities are
excluded by Theorem~\ref{thm:exceptional}.

Finally, suppose that $3\leq c_2<q$, where $q=k-1$.
Recall that $A_u$ acts faithfully on $\Gamma(u)$ in its
natural projective-line action, with $\PSL(2,q)\leq A_u\leq\PGammaL(2,q)$.
Thus the classification in~\cite[p.~4, Case~(IV)]{CP1983} applies.
The parameters $v$ and $k$ in that
reference correspond to our $k$ and $c_2$, respectively.
It follows that
\begin{enumerate}[\rm (i)]
\item
$(k,c_2)=(6,3)$ and $\Gamma\cong B'(H_{11})$, as part (5);
\item
$(k,c_2)=(8,4)$ and $\Gamma\cong B'(PG(3,2))$.
\end{enumerate}
By Lemma~\ref{exam:B32}, case~(ii) cannot occur
because $B'(PG(3,2))$ is not $2$-distance-primitive.

Conversely, Proposition~\ref{examples} and
Lemmas~\ref{exam:Hamming}--\ref{exam:B11} show that all the graphs
listed in Theorem~\ref{thm:main} are $2$-distance-primitive.
This completes the proof. \qed

\section*{Acknowledgements}
The second author was partially supported by the National Natural Science Foundation of China (No. 12401463),
and the third  author was partially supported by the National Natural Science Foundation of China (Nos. 12501469 and 12671398).

\end{document}